\documentclass[a4paper,10pt]{amsart}

\usepackage{amsmath,amssymb,amsthm,a4wide}

\usepackage{tikz}
\usetikzlibrary{shapes}
\usetikzlibrary{intersections}

\usepackage{graphicx}
\usepackage{graphics}

\newtheorem{theorem}{Theorem}

\begin{document}

\title[A spectral gap estimate and applications]{A spectral gap estimate and applications}

\author[]{Bogdan Georgiev}
\address{Max Planck Institute for Mathematics, Vivatsgasse 7, 53111 Bonn, Germany}
\email{bogeor@mpim-bonn.mpg.de}

\author[]{Mayukh Mukherjee}
\address{Max Planck Institute for Mathematics, Vivatsgasse 7, 53111 Bonn, Germany}
\email{mukherjee@mpim-bonn.mpg.de}

\author[]{Stefan Steinerberger}
\address{Department of Mathematics, Yale University, 06511 New Haven, CT, USA}
\email{stefan.steinerberger@yale.edu}

\begin{abstract} We consider the Schr\"odinger operator
 $$-\frac{d^2}{d x^2} + V \qquad \mbox{on an interval}~~[a,b]~\mbox{with Dirichlet boundary conditions},$$
where $V$ is bounded from below and prove a lower bound on the first eigenvalue $\lambda_1$ in terms of sublevel estimates: if
 $ w_V(y) = |I_y|,\text{  where  }  I_y := \left\{  x \in [a,b]:  V(x) \leq y \right\},$
	then
$$ \lambda_1 \geq    \frac{1}{250} \min_{y > \min V}{\left(\frac{1}{w_V(y)^2} + y\right)}.$$
The result is sharp up to a universal constant if $\left\{  x \in [a,b]:  V(x) \leq y \right\}$ is an interval for the value of $y$ solving the minimization problem.
An immediate application is as follows: let $\Omega \subset \mathbb{R}^2$ be a convex domain with inradius $\rho$ and diameter $D$ and let $u:\Omega \rightarrow \mathbb{R}$ be the first eigenfunction of the Laplacian $-\Delta$ on $\Omega$ with Dirichlet boundary conditions on $\partial \Omega$. We prove
$$ \| u \|_{L^{\infty}} \lesssim \frac{1}{\rho^{}} \left( \frac{\rho}{D} \right)^{1/6} \|u\|_{L^2},$$
which answers a question of van den Berg in the special case of two dimensions.
\end{abstract}


\maketitle

\section{Introduction and statement of results}

\subsection{Introduction and Result}
We consider one-dimensional Schr\"odinger operators of the form
 $$-\frac{d^2}{d x^2} + V \qquad \mbox{on an interval}~~[a,b]~\mbox{with Dirichlet boundary conditions},$$
 where $V$ is assumed to be bounded from below. It is intuitively clear that the smallest eigenvalue $\lambda_1$ will mainly depend on the minimal value attained by the potential as value as the growth rate around that minimal value. We
make this intuition precise. In the case where all sublevel sets are intervals (i.e. which, for instance, is the case for convex potentials $V$), this
lower bound is sharp up to a universal multiplicative constant.
 We define a function $w:\{y \in \mathbb{R} : y \geq \min V\} \rightarrow \mathbb{R}^+$  by measuring the length of sublevel sets via
 $$ w_V(y) = |I_y|,\text{  where  }  I_y := \left\{  x \in [a,b]:  V(x) \leq y \right\}.$$
 We will prove that a lower bound on the first eigenvalue can be given as the solution of a minimization problem involving the function $w(y)$. 

 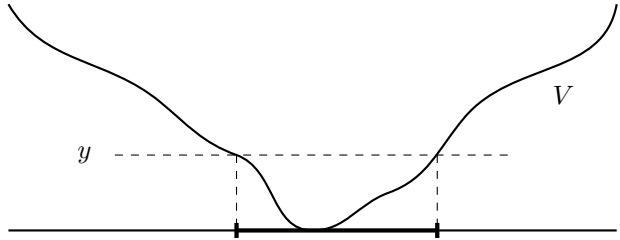
\begin{figure}[h!]
 	\begin{center}
 		\begin{tikzpicture}[scale=1]
 		\draw [thick] (-4,0) -- (4,0);
 		\draw [thick] (-4,3) to[out=300, in=140] (-2,1.7)  to[out=320, in=160] (-1,1) to[out=340, in=180] (0,0)  to[out=0, in=200] (1,0.5)    to[out=20, in=230] (2,1.5)    to[out=50, in=260] (4,3) ; 
 		\node at (-3, 1) {$y$};
 		\node at (3.3, 1.8) {$V$};
 		\draw [dashed] (-2.6, 1) -- (2.6, 1);
 		\draw [dashed] (-1, 1) -- (-1,0);
\draw [ultra thick] (-1, 0) -- (1.64,0);
 		\draw [dashed] (1.64, 1) -- (1.64,0);
 		\draw [ultra thick] (1.64, -0.1) -- (1.64,0.1);
 		\draw [ultra thick] (-1, -0.1) -- (-1,0.1);
 		\end{tikzpicture}
 	\end{center}
 	\caption{$w_V(y)$ is the length of the interval $\left\{x: V(x) \leq y\right\}$.}
 \end{figure}

 \begin{theorem} \label{th:Eignval-Estimate} If $V$ is bounded from below, then the smallest eigenvalue $\lambda_1$ of $-d^2/dx^2 + V$ with Dirichlet conditions at the endpoints of the interval satisfies
 	$$ \lambda_1 \geq    \frac{1}{250} \min_{y > \min V}{\left(\frac{1}{w_V(y)^2} + y\right)}.$$
 	Moreover, if the sublevel set $\left\{x: V(x) \leq y^*\right\}$ is an interval for the value $y^* \in \mathbb{R}$ solving the minimization problem, then
 	$$ \lambda_1 \leq \pi^2 \min_{y > \min V}{\left(\frac{1}{w_V(y)^2} + y\right)}.$$
 \end{theorem}
 The upper bound is quite simple and follows from a testing argument that actually implies the
 slightly sharper result
 $$\lambda_1\leq \min_{y > \min V}{\left(\frac{ \pi^2}{w_V(y)^2} + y\right)}.$$
 The constant $1/250$ in the lower bound seems far from optimal and it could be of interest to obtain better results. 
 It could also be of interest to study higher-dimensional analogues.
While we do not know of any such result in the literature, a special case of the result has been established in work of Jerison \cite{jer3} on the
ground state of the Laplacian in convex domains in $\mathbb{R}^2$: the eigenvalue (and the scale of localization) is determined in terms of
a geometric characterization that is equivalent to Theorem 1, details are given at the end of the paper.

\subsection{An application}
Let $\Omega \subset \mathbb{R}^n$ be a convex domain with inradius $\rho$ and diameter $D$ and let $u \in L^2(\Omega)$ be
the ground state of the Laplacian $-\Delta$ in $\Omega$ with Dirichlet boundary conditions on $\partial \Omega$. A classical inequality
of Chiti \cite{chi1, chi2} states that
$$
 \| u \|_{L^{\infty}} \lesssim \frac{1}{\rho^{n/2}}   \| u \|_{L^{2}},
$$
where the sharp constant is known and assumed for the ball. A natural question, asked by van den Berg \cite{vdb}, is whether this inequality
could be improved when the diameter $D$ gets large -- motivated by explicit computations on cones, he conjectured that
$$ \| u \|_{L^{\infty}} \lesssim \frac{1}{\rho^{n/2}} \left( \frac{\rho}{D} \right)^{1/6} \| u \|_{L^2}.$$

The second main result of this paper is to prove the conjecture in two dimensions.

\begin{theorem} \label{th:van-den-Berg-Conj} Let $\Omega \subset \mathbb{R}^2$ be convex with diameter $D$ and inradius $\rho$ and let $u:\Omega \rightarrow \mathbb{R}$ be the ground state
of the Laplacian with Dirichlet boundary conditions. Then
$$
	 \| u \|_{L^{\infty}} \lesssim \frac{1}{\rho^{}} \left( \frac{\rho}{D} \right)^{1/6} \|u \|_{L^2}.
$$
\end{theorem}

Roughly, our proof consists of two ingredients: we use a result of Grieser \& Jerison \cite{jer1}
to reduce the problem to the problem of analyzing the ground state of a Schr\"odinger operator $-d^2/dx^2 + V$ on a compact interval;
further, we use Theorem \ref{th:Eignval-Estimate} to estimate the first eigenvalue and first eigenfunction of such an operator.

\section{Proofs}

\subsection{Proof of Theorem \ref{th:Eignval-Estimate}}
\begin{proof} Assuming that 
	$$  \left\{x: V(x) \leq y\right\} \qquad \mbox{is an interval},$$
	the upper bound is fairly easy and follows immediately from choosing a using suitable test function in the Rayleigh-Ritz quotient: for every $y > \min V$, define
	$$ f(x) = \begin{cases}  \sin{ \left(  \pi \frac{x - \inf I_{y}}{|I_{y}|}  \right)} \qquad &\mbox{if}~x \in I_y = \left\{x: V(x) \leq y\right\} \\
	0 \qquad &\mbox{otherwise.} \end{cases}$$
	A direct computation using the Rayleigh-Ritz quotient shows
	\begin{align*}
	\lambda_1  \int_{a}^{b}{ f(x)^2 dx} &\leq \int_{a}^{b}{ f'(x)^2 + V(x) f(x)^2 dx} \\
	&\leq \frac{\pi^2}{w(y)^2}  \int_{a}^{b}{ f(x)^2 dx} +  y \int_{a}^{b}{ f(x)^2 dx} \leq \pi^2 \left( \frac{1}{w(y)^2} + y\right)\int_{a}^{b}{ f(x)^2 dx}.
	\end{align*}
	Since this is true for all $y$, we can take the minimum of the arising expression as an upper bound on the first eigenvalue, which implies the upper bound on the eigenvalue.\\

	Let now $f^*$ denote the symmetrically
	decreasing rearrangement of $f$ around $(a+b)/2$ and let $V_{*}$ denote the symmetrically increasing rearrangement of $V$ around $(a+b)/2$. The Hardy-Littlewood rearrangement
	inequality implies
	$$ \int_{a}^{b}{ V(x) f(x)^2 dx} \geq  \int_{a}^{b}{ V_*(x) f^*(x)^2 dx}$$
	and the Polya-Szeg\H{o} inequality implies
	$$  \int_{a}^{b}{f(x)^2 dx} = \int_{a}^{b}{f^*(x)^2 dx}  \qquad \mbox{and} \qquad  \int_{a}^{b}{ \left(\frac{d}{dx} f^*(x) \right)^2 dx} \leq  \int_{a}^{b}{ \left(\frac{d}{dx} f(x) \right)^2 dx}.$$
	From this, we can infer that
	$$ \lambda_1\left( - \frac{d^2}{dx^2} + V^* \right) \leq \lambda_1\left( - \frac{d^2}{dx^2} + V \right)$$
	and will now proceed to bound the smallest eigenvalue of the rearranged potential from below.
	The argument will keep track of three a priori unspecified
	real constants $\alpha, \beta, \gamma$ (an optimization over which will then yield an explicit lower bound). 
	Let $f$ denote the $L^2-$normalized ground state of $-\Delta + V$ with Dirichlet conditions at $a$ and $b$. 
	It is classical (and follows from the ordinary differential equation and the unimodal structure of $V$) that the ground state is first monotonically increasing and then monotonically decreasing.
	Let $J \subset [a,b]$ be the smallest interval such that
	$$ \int_{J}{f(x)^2 dx} \geq \alpha \|f\|^2_{L^2},$$
	where $\alpha \in (0, 1)$ is a constant that will be determined later. We start by observing that the function $f$ 
	has the same value on both boundary points of $J$: if it did not, then continuity of $f$ suffices to ensure that we could slightly
	slide the interval to increase the contained $L^2-$norm, which
	would then allow shrinking the interval, thus contradicting the minimality of $J$.
	
	We consider the interval 
	$$J^* = \left[\frac{a+b}{2} - \frac{|J|}{2}, \frac{a+b}{2} + \frac{|J|}{2} \right] \quad \mbox{and observe that} \quad
	\int_J f(x)^2 dx = \int_{J^*} f^*(x)^2 dx = \alpha \|f \|_{L^2}^2.
	$$
	Clearly, the $L^2-$normalization of $f$ implies that
	$$ \int_{[a,b] \setminus J^*}^{}{ f^*(x)^2 dx} = (1 - \alpha) \|f\|^2_{L^2}$$
	and therefore
	\begin{align*}
	\int_{a}^{b}{ V_*(x) f^*(x)^2 dx} &\geq \int_{[a,b]\setminus J^*}{ V_*(x) f^*(x)^2 dx}   \\
	&\geq  w^{-1}_*(|J^*|) (1- \alpha)  \|f\|^2_{L^2} =  (1- \alpha) w^{-1}(|J|) \|f\|^2_{L^2}.
	\end{align*}
	The remainder of the proof deals with the gradient term and proceeds via case distinction: either there is a lot of oscillation on the interval $J$ or there is not. More
	precisely, for a constant $\beta > 0$ to be determined later, we are either dealing with
	$$ \mbox{either}\quad \int_{J}{f'(x)^2 dx} \geq   \frac{\beta}{|J|^2}  \int_{J}{f(x)^2 dx}  \qquad \mbox{or} \qquad \int_{J}{f'(x)^2 dx} \leq   \frac{\beta}{|J|^2}  \int_{J}{f(x)^2 dx}.$$
	
	\textbf{Case 1} (Lots of Oscillation)\textbf{.} As it turns out, this case is easy to deal with since we can estimate
	$$ \int_{a}^{b}{ f'(x)^2 dx} \geq \int_J { f'(x)^2 dx} \geq  \frac{\beta}{|J|^2}  \int_{J}{f(x)^2 dx} =  \frac{\alpha \beta}{|J|^2} \|f\|_{L^2}^2.$$
	Altogether, this means that we have
	\begin{align*}
	\int_{a}^{b}{ f'(x)^2 + V(x) f(x)^2 dx} &\geq \left( \frac{\alpha \beta}{|J|^2} +  (1-\alpha) w^{-1}(|J|) \right)\|f\|_{L^2}^2\\
	&\geq \min(\alpha \beta, 1 - \alpha) \left( \frac{1}{|J|^2}  + w^{-1}(|J|) \right)\|f\|_{L^2}^2.
	\end{align*}
	However, a change of variables immediately shows that
	$$  \left( \frac{1}{|J|^2}  + w^{-1}(|J|) \right) \geq  \min_{y > \min V}{\left(\frac{1}{w(y)^2} + y\right)},$$
	which is the desired statement with a constant $\min(\alpha \beta, 1 - \alpha)$.\\
	
	\textbf{Case 2} (Almost flat)\textbf{.} The remaining case is where $f$ has relatively little oscillation on $J$ 
	$$ \int_{J}{f(x)^2 dx} = \alpha \|f\|^2_{L^2} \qquad \mbox{and} \qquad \int_{J}{f'(x)^2 dx} \leq \frac{\beta}{|J|^2}  \int_{J}{f(x)^2 dx}.$$
	The remainder of the argument is as follows: we will show that the function has to decay outside of $J$ with at least a certain speed
	(otherwise, if the function were to remain close to constant, it would eventually violate the $L^2-$normalization of $f$).
	As proved above, $f$ has the same value at the two endpoints of $J$, which we will denote by $\varepsilon := f(\partial J)$. We start by showing that this
	number cannot be too small and consider $g = f - \varepsilon$ on $J$. Using the Cauchy-Schwarz inequality, we
	see that
	\begin{align*}
	\int_{J}{g(x)^2 dx} &= \int_{J}{(f(x) - \varepsilon)^2 dx} \geq \int_J f(x)^2 dx - 2 \varepsilon \int_{J}{f(x)dx} \\
	&\geq \alpha - 2 \varepsilon |J|^{1/2} \left( \int_{J}{f(x)^2 dx} \right)^{1/2}\\
	&\geq \alpha - 2 \varepsilon |J|^{1/2} \alpha^{1/2}.
	\end{align*}
	At the same time, we have that
	$$ \int_{J}{g'(x)^2 dx} = \int_{J}{f'(x)^2 dx} \leq  \frac{\beta}{|J|^2}  \int_{J}{f(x)^2 dx} = \frac{\alpha \beta}{|J|^2}.$$
	The classical Poincar\'{e} inequality for the Dirichlet-Laplacian on an interval implies that
	$$  \int_{J}{g'(x)^2 dx}  \geq \frac{\pi^2}{|J|^2}   \int_{J}{g(x)^2 dx}  $$
	and therefore
	$$ \frac{\alpha \beta}{|J|^2} \geq \frac{\pi^2}{|J|^2}  (\alpha - 2\varepsilon|J|^{1/2} \alpha^{1/2} )$$
	and after rearrangement
	$$ \varepsilon \geq  \frac{\alpha^{1/2}}{2|J|^{1/2}}  \left(1- \frac{\beta}{\pi^2}\right).$$
	We now show that $f$ has to decay at least with a certain speed outside of $J$, otherwise there would be too much $L^2-$mass outside of $J$. 	
	We fix another positive constant $\gamma$ and look at the values of $f$ at points, whose distance to $J$ is $\gamma |J|$. Assuming that $J = [j_1, j_2]$ and using
	the monotonicity of $f$, we have
	\begin{align*}
	1 = \|f\|^2_{L^2([a,b])} \geq \int_{j_1}^{j_2+\gamma |J|}{ f(x)^2 dx} \geq   (1 + \gamma)|J| f(j_2 + \gamma |J|)^2
	\end{align*}
	and thus $$ f(x) \leq \frac{1}{(1+\gamma)^{1/2}|J|^{1/2}} \qquad \mbox{for all}~x~\mbox{at distance at least}~\gamma|J|~\mbox{from the interval}~J.$$
	
	Using Euler-Lagrange equations it is a basic exercise to show
	$$ \inf \left\{ \int_{c}^{d}{h'(x)^2 dx} \big|~h(c)=C, h(d) = D \right\} = \frac{(D-C)^2}{d-c}$$
	because the minimizer is simply the linear function with appropriate values at the endpoints.

	For our problem, this implies that
	\begin{align*}
	\int_{a}^{b}{f'(x)^2 dx} &\geq \int_{j_2}^{j_2 + \gamma |J|} f'(x)^2 dx \\
	&\geq \int_{j_2}^{j_2 + \gamma |J|} h'(x)^2 dx \geq \frac{1}{\gamma |J|} \left( \varepsilon - f(j_2 + \gamma |J|) \right)^2 \\ 
	&\geq \frac{1}{\gamma |J|^2} \left[ \frac{\alpha^{1/2}}{2}  \left(1- \frac{\beta}{\pi^2}\right) - \frac{1}{(1+\gamma)^{1/2}} \right]^2 =: \frac{1}{|J|^2} F(\alpha, \beta, \gamma),
	\end{align*}
	where we have to assume that
	\begin{align*}
	\frac{\alpha^{1/2}}{2}  \left(1- \frac{\beta}{\pi^2}\right) \geq \frac{1}{(1+\gamma)^{1/2}}.
	\end{align*}
	
	We can argue as in Case 1 and obtain that
	\begin{align*}
	\int_{a}^{b}{ f'(x)^2 + V(x) f(x)^2 dx} &\geq \left( \frac{1}{|J|^2} F(\alpha, \beta, \gamma) +  (1 - \alpha) w^{-1}(|J|) \right)\|f\|_{L^2}^2\\
	&\geq \min \left( F(\alpha, \beta, \gamma \right), (1-\alpha)) \left( \frac{1}{|J|^2}  + w^{-1}(|J|) \right)\|f\|_{L^2}^2,
	\end{align*}
	which is the desired result.\\
	
	\textbf{The numerical constant.} We conclude by untangling the relationship between the implicit constants.
	Altogether, we have the following system of constraints 
	\begin{align*} \label{eq:Constraints}
	\alpha \in (0, 1), \quad \beta \in (0, \pi^2), \quad \gamma \geq \frac{4}{\alpha} \left(1 - \frac{\beta}{\pi^2} \right)^{-2} - 1.
	\end{align*}
	and are trying to estimate
	$$ \max_{\alpha, \beta, \gamma}  \min \left\{ \frac{1}{\gamma}\left[ \frac{\alpha^{1/2}}{2}  \left(1- \frac{\beta}{\pi^2}\right) - \frac{1}{(1+\gamma)^{1/2}} \right]^2, 1-\alpha, \alpha \beta \right\}.$$
	The lower bound $1/250$ follows from setting
	$$ \alpha = \frac{99}{100}, ~\beta= \frac{7}{1000} \quad \mbox{and} \quad \gamma = 14.1327.$$
	
	\textbf{$L^{\infty}-$bounds.} We finally establish bounds on the $L^{\infty}-$norm of the eigenfunction in the case where the potential is nonnegative $V \geq 0$,
	which is an easy combination of the Cauchy-Schwarz inequality and the fact that the ground state
	vanishes on the endpoints. More precisely, for every $a \leq x \leq b$ we distinguish between the cases
	$$ \int_{a}^{x}{f(y)^2 dy} \leq \frac{1}{2} \qquad \mbox{or} \qquad  \int_{x}^{b}{f(y)^2 dy} \leq \frac{1}{2}.$$
	In the first case, we estimate
	\begin{align*}
	f(x)^2 = f(x)^2 - f(a)^2 &= \int_{a}^{x}{ \left(\frac{d}{dy} f(y)^2\right) dx} = \int_{a}^{x}{ 2 f(y) f'(y) dy} \\
	&\leq 2 \left( \int_{a}^{x}{f(y)^2 dy}\right)^{1/2}  \left( \int_{a}^{x}{f'(y)^2 dy}\right)^{1/2} \\
	&\leq \sqrt{2}  \left( \int_{a}^{b}{f'(x)^2 + V(x) f(x)^2 dx}\right)^{1/2} \leq \sqrt{2\lambda_1}.
	\end{align*}
	In the second case, we change signs and argue
	\begin{align*}
	f(x)^2 = f(x)^2 - f(b)^2 &= -\int_{x}^{b}{ \left(\frac{d}{dy} f(y)^2\right) dx} = -\int_{x}^{b}{ 2 f(y) f'(y) dy} \\
	&\leq 2 \left( \int_{x}^{b}{f(y)^2 dy}\right)^{1/2}  \left( \int_{x}^{b}{f'(y)^2 dy}\right)^{1/2} \\
	&\leq \sqrt{2}  \left( \int_{a}^{b}{f'(x)^2 + V(x) f(x)^2 dx}\right)^{1/2} \leq \sqrt{2\lambda_1}
	\end{align*}
	and therefore
	$$
	\|f\|_{L^{\infty}} \leq (2 \lambda_1)^{1/4}.
	$$
\end{proof}

\subsection{Proof of Theorem \ref{th:van-den-Berg-Conj}}
\begin{proof} We begin by giving a short overview of the overall idea. We can use scaling to restrict ourselves to convex domains with inradius $\rho = 1$. 
	We first recall 
	Theorem 1.6 of \cite{jer1}: assume we rotate the convex domain $\Omega$ so that the projection onto the $y$-axis has least length and dilate (if necessary) so that this length is $1$. The boundary of $\Omega$ can then be written as the union of the graphs of two functions $f_1(x) \leq f_2(x)$ on $[a,b]$, which satisfy 
	$$0 \leq f_1(x) \leq f_2(x) \leq 1 \text{  for }  a \leq x \leq b,$$
	$$ \min_{a \leq x \leq b}f_1(x) = 0, \quad \max_{a \leq x \leq b}f_2(x) = 1.$$
	Their result then states that the profile of the eigenfunction $u$ on $\Omega$ is essentially
given by the ground state $\phi$ of the Schr\"odinger operator
$$ -\frac{d^2}{dx^2} + \frac{\pi^2}{h(x)^2}, \quad \mbox{where}~h(x) := f_2(x) - f_1(x).$$ 
Formally, let $L$ be the length of the longest interval $I \subset [a, b]$ such that
$$
	h(x) \geq 1 - \frac{1}{L^2},
$$
on $I$ and let
$$
	\alpha(x,y) = \pi \frac{y - f_1(x)}{h(x)}.
$$
\begin{quote}
\begin{theorem}[\cite{jer1}, Theorem 1.6]\label{thm:GJ}
	Normalize $u$ and $\phi$, so that $\max u = \max \phi = 1$. Then there is an absolute constant $C$ such that
$$
		| u(x,y) - \phi(x) \sin \alpha(x,y)| \leq \frac{C}{L},
$$
	for all $x \in I'$ where $I'$ is the interval concentric with $I$ of half the length.
\end{theorem}
\end{quote}

Our proof will now proceed as follows: we first study exclusively the ground states of Schr\"odinger operators that
can possibly arise from convex domains and establish some basic properties. Afterwards, we use the 
above result to transfer the results to the profile of the eigenfunction on the convex domain.\\

\textit{Step 1} (Bound on the ground state). \\
We make use of Theorem 2 and replace our study of the first eigenvalue of the Schr\"odinger operator associated to the domain with the study of
$$\min_{y > \min V}{\left(\frac{1}{w_V(y)^2} + y\right)}.$$
It is easy to see that a larger potential in the Schr\"odinger operator $V_1 \geq V_2$ leads to shorter sublevels $w_{V_1}(y) \leq w_{V_2}(y)$, which can only increase the
minimal value in the minimization problem. Every convex body of inradius 1 and diameter $D$ contains a disk of radius $1$ and a point at distance $\sim D/2$,
which gives rise to a cone by convexity -- this cone is universal in the sense that it is strictly contained in every other convex domain with inradius 1
and diameter $D$ and its potential therefore dominates all other potentials. 
\begin{center}
	\begin{figure}[h!]
		\begin{tikzpicture}
		\draw [thick] (0,0) circle (1cm);
		\draw [thick] (0,1)--(6,0);
		\draw [thick] (0,-1)--(6,0);
		\end{tikzpicture}
		\caption{A description of the object contained in every convex set with inradius 1 and diameter $D$: a ball with inradius $1 $ and a cone with height $(D-1)/2$.}
	\end{figure}
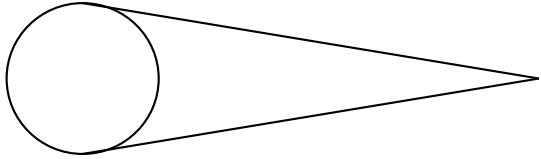
\end{center}
A simple computation shows that in this case the relevant potential is, up to universal constants,
$$ V \sim \frac{D^2}{(D-x)^2} \qquad \mbox{on}~[0,D].$$
We now use the fact that adding or subtracting constants to the potential does not
change the ground state and replace the potential under consideration by
$$ V \sim \frac{D^2}{(D-x)^2} - 1$$
and apply Theorem 2 to this potential instead.
 We then obtain
 $$ w_{V}(y) = D\left(1-\frac{1}{\sqrt{1 + y}}\right)$$
 and thus
$$\left(\frac{1}{w_{V}(y)^2} + y\right)  \sim D^{-2}\left(1 - \frac{1}{\sqrt{1 + y}}\right)^{-2} + y \sim \frac{1}{D^2 y^2} + y$$
which is minimized for $y \sim D^{-2/3}$ and gives the upper bound $\lambda_1 \lesssim D^{-2/3}$. A classical $L^{\infty}-$bound, proven
for convenience of the reader at the end of the paper, gives
$$ \| \phi\|_{L^{\infty}} \lesssim \lambda_1^{1/4} \| \phi \|_{L^2}.$$
This gives the desired result $\| \phi\|_{L^{\infty}} \lesssim D^{-1/6} \| \phi \|_{L^2}$ for the `profile'-eigenfunction of the Schr\"odinger operator.
The second part of Theorem 2 implies the following: if we take $J$ to be the shortest interval containing half of the $L^2-$mass of $\phi$, then
$$ |J| \gtrsim \frac{1}{\sqrt{\lambda_1}}.$$
These two facts imply that $\phi$ is essentially constant on an interval of length $L \sim 1/\sqrt{\lambda_1}$ and that a constant proportion of $L^2-$mass is contained on that interval.\\

\textit{Step 2} (Using Grieser-Jerison). 
A simple computation (mentioned in \cite{jer1}, see also \cite[Lemma 2.4]{jer3}) shows that 
the comparison is accurate at least on length scale $\sim \lambda_1^{-1/2}$ around the maximum. Suppose now that the functions $u_{1}, \phi_1$ are rescaled 
from the eigenfunctions in such a way that
$$\max u_1 = \max \phi_{1} = 1$$
and that
$$| u_1(x,y) - \phi_1(x) \sin \alpha(x,y)| \leq \frac{C}{L}.$$
We know from the first part that $\phi$ is essentially constant on scale $\gtrsim \lambda_1^{-1/2}$, that $\| \phi \|_{L^{\infty}} \lesssim \lambda_1^{1/4}$ and that a
constant proportion of the $L^2-$norm is there. This implies for the
rescaled function $\phi_1$ that
$$ \int_{J}{ \phi_1(x)^2 dx} \gtrsim  \int_{J}{ \lambda_1^{-1/2} \phi(x)^2 dx} \sim \lambda_1^{-1/2},$$
where $J$ is the interval of length $\lambda_{1}^{-1/2}$ on which $\phi$ is essentially constant. 
This, however, implies 
$$ \int_{\Omega}{u_1(x,y)^2 dx dy} \gtrsim \lambda_1^{-1/2}$$
and therefore
$$  \|u_1\|^2_{L^{\infty}} = 1  \lesssim \lambda_{1}^{1/2} \|u_1\|_{L^2}^2 \lesssim  D^{-1/3} \|u_1\|_{L^2}^2$$
which is the desired result.
\end{proof}

We emphasize a useful connection between the estimate
$$ \lambda_1\sim  \min_{y > \min V}{\left(\frac{1}{w_V(y)^2} + y\right)}$$
and the length $L$ of the longest interval $I \subset [a, b]$ such that
$$h(x) \geq 1 - \frac{1}{L^2}$$
in the Grieser-Jerison Schr\"odinger operator \cite{jer1, jer3}. Note that
$$ L \leq \left|\left\{x:  \frac{1}{h(x)^2} - 1 \leq \frac{1}{L^2}\right\}\right| \qquad \mbox{and thus} \qquad w_V(L^{-2}) \sim  L.$$
As a consequence, for $y = L^{-2}$, we observe that
$$\frac{1}{w_V(y)^2} = \frac{1}{w_V(L^{-2})^2} \sim L^{-2} = y$$
which corresponds to a balancing of terms in the functional and immediately yields that $\lambda_1 \sim L^{-2}$.

\subsection*{Acknowledgements} The first and second authors gratefully acknowledge the Max Planck Institute for Mathematics, Bonn for providing ideal working conditions.
The third author was partially supported by an AMS Simons Travel grant and INET Grant \#INO15-00038

\end{document}